\documentclass[reqno]{conm-p-l}
\usepackage{amssymb,amsmath,amscd,euscript,mathrsfs}

\usepackage{graphicx}
\usepackage{amsfonts}

\usepackage[all]{xy}

\input epsf

\theoremstyle{plain}
\numberwithin{equation}{section}
\newtheorem{theorem}{Theorem}[section]

\theoremstyle{definition}
\newtheorem{definition}[theorem]{Definition}

\newtheorem*{ack}{Acknowledgements}

\def\fa{{\mathfrak a}}

\def\fg{{\mathfrak g}}
\def\fh{{\mathfrak h}}

\def\fk{{\mathfrak k}}

\def\fs{{\mathfrak s}}

\newcommand{\bc}{\mathbf{c}}

\newcommand{\rO}{\mathrm{O}}
\renewcommand{\Re}{\mathrm{Re}\,}
\def\Im{{\Rm Im}\,}

\def\cA{{\mathcal A}}

\def\cE{{\mathcal E}}
\def\cF{{\mathcal F}}
\def\cH{{\mathcal H}}

\def\cP{{\mathcal P}}
\def\cR{{\mathcal R}}
\def\cS{{\mathcal S}}

\newcommand{\C}{\mathbb{C}}
\newcommand{\D}{\mathbb{D}}
\newcommand{\T}{\mathbb{T}}
\newcommand{\R}{\mathbb{R}}

\newcommand{\Z}{\mathbb{Z}}

\newcommand{\rS}{\mathrm{S}}
\newcommand{\rU}{\mathrm{U}}

\newcommand{\wG}{\widehat{G}_K}

\newcommand{\wf}{\widehat{f}}
\newcommand{\wtf}{\widetilde{f}}

\newcommand{\wphi}{\widehat{\phi}}
\newcommand{\Tr}{\mathrm{Tr}}

\newcommand{\SU}{\mathrm{SU}}
\newcommand{\SO}{\mathrm{SO}}

\newcommand{\hm}{\widehat{\mu}}

\newcommand{\SL}{\mathrm{SL}}
\newcommand{\ip}[2]{\langle#1,#2 \rangle}

\def\sideremark#1{\ifvmode\leavevmode\fi\vadjust{\vbox to0pt{\vss
 \hbox to 0pt{\hskip\hsize\hskip1em%
 \vbox{\hsize2cm\tiny\raggedright\pretolerance10000
 \noindent #1\hfill}\hss}\vbox to8pt{\vfil}\vss}}}

\newcommand{\rank}{\mathrm{rank}}
 
\newcommand{\rE}{\mathrm{E}}
\newcommand{\rG}{\mathrm{G}}
\newcommand{\rH}{\mathrm{H}}
\newcommand{\vv}{\varphi^\vee}
\newcommand{\vx}{x^\vee}
\renewcommand{\Im}{\mathrm{Im}}

\usepackage[icelandic,english]{babel}
\usepackage[T1]{fontenc}

\begin{document}

\title{On the Life and Work of S. Helgason}

\author{G. \'Olafsson}
\address{Department of Mathematics, Louisiana State University, Baton Rouge,
LA, 70803 USA} \email{olafsson@math.lsu.edu}

\author{R. J. Stanton}
\address{Department of Mathematics, Ohio State University, Columbus, OH, 43210 USA} \email{stanton@math.ohio-state.edu}

\thanks{G.\'Olafsson acknowledges the support of NSF Grant DMS-1101337 during the preparation of this article.}


\date{February 26, 2013}

\begin{abstract}
This article is a contribution to a Festschrift for S. Helgason. After a biographical sketch, we survey some of his research on several topics in geometric and harmonic analysis during his long and influential career. While not an exhaustive presentation of all facets of his research, for those topics covered we include reference to the current status of these areas. 
\end{abstract}

\maketitle
\noindent\section*{Preface}
\noindent

Sigur{\dh}ur Helgason is known worldwide for his first book \textit{Differential Geometry and Symmetric Spaces}. With this book he provided an entrance to the opus of \'Elie Cartan and Harish-Chandra to generations of mathematicians. On this the occasion of his 85th birthday we choose to reflect on the impact of Sigur{\dh}ur Helgason's  sixty years of mathematical research. He was among the first to investigate systematically the analysis of differential operators on reductive homogeneous spaces. His research on Radon-like transforms for homogeneous spaces presaged the resurgence of activity on this topic and continues to this day. Likewise he gave a geometrically motivated approach to harmonic analysis of symmetric spaces. Of course there is much more -  eigenfunctions of invariant differential operators, propagation properties of differential operators, differential geometry of homogeneous spaces, historical profiles of mathematicians. Here we shall present a survey of some of these contributions, but first a brief look at the man.

\section{Short Biography}
Sigur{\dh}ur Helgason was born on September 30, 1927 in Akureyri, in northern Iceland. His parents were Helgi Sk\'ulason (1892-1983)  and Kara Sigur{\dh}ard{\'o}ttir Briem (1900-1982), and he had a brother Sk\'uli Helgason (1926-1973) and a sister Sigri{\dh}ur Helgad{\'o}ttir (1933-2003). 
Akureyri was then the second largest city in Iceland with about 3,000 people living there, whereas the population of Iceland was about 103,000. As with other cities in northern Iceland, Akureyri was isolated, having only a few roads so that horses or boats were the transportation of choice. Its schools, based on Danish traditions, were good. The Gymnasium in
Akureyri was established in 1930 and was the second Gymnasium in Iceland. There Helgason studied mathematics,
physics, languages, amongst other subjects during  the years 1939-1945.  He then went to the University of
Iceland  in Reykjav{\'i}k where he enrolled in the school of engineering, at that time the only way there to study mathematics. In 1946 he began studies at the
University of  
Copenhagen from which he received the Gold Medal in 1951 for his  work on Nevanlinna-type value distribution theory for analytic almost-periodic functions. His paper on the subject became his master's thesis in 1952. Much later a summary appeared in \cite{He89}.

Leaving Denmark in 1952  he went to Princeton University to complete his graduate studies. He received a  Ph.D.  in 1954 with the thesis, \textit{Banach Algebras and Almost Periodic Functions},
under the supervision of Salomon Bochner.

He began his professional career as a C.L.E. Moore Instructor at M.I.T. 1954-56. After leaving Princeton his interests had started to move towards
two areas that remain the main focus of his research.
The first, inspired by Harish-Chandra's ground breaking work on the
representation theory of semisimple Lie groups, was Lie groups and harmonic analysis on
symmetric spaces;  the second was  the Radon transform, the motivation having come from
reading the page proofs of Fritz
John's famous 1955  book {\it Plane Waves and Spherical Means}. He returned to
Princeton for 1956-57  where his interest in Lie groups and
symmetric spaces  led to his first work on applications of Lie theory to
differential equations, \cite{He59}. He moved to the University of Chicago for 1957-59, where he started
work on his first book \cite{He62}.  He then went  to
Columbia University for the fruitful period 1959-60, where he shared an office with Harish-Chandra. In 1959 he  joined the faculty at M.I.T.
where he has remained these many years, being full professor since 1965. The periods 
1964-66, 1974-75, 1983 (fall) and 1998 (spring) he spent at the 
Institute for Advanced Study, Princeton, and the periods 1970-71
and 1995 (fall) at the Mittag-Leffler Institute, Stockholm. 

He has 
been awarded a degree Doctoris Honoris Causa by several universities, notably the University of Iceland, the University of
Copenhagen and the University of Uppsala. In 1988 the American Mathematical Society awarded him the
Steele Prize for expository writing citing his book \textit{Differential Geometry and Symmetric Spaces} and its sequel.
Since 1991 he carries the Major Knights Cross of the Icelandic Falcon.

\section{Mathematical Research}

In the Introduction to his selected works, \cite{Sel},  Helgason himself gave a personal description of his work and
how it relates to his published articles.  We recommend this for the clarity of exposition we have come to expect from him as well
as the insight it provides to his motivation. An interesting interview with him also may be found in \cite{Seg}.
Here we will  discuss parts of this work, mostly those familiar to us. We start with his work on invariant differential operators, continuing with his work on Radon transforms, his work related to symmetric spaces and representation theory, then a sketch of his work on wave equations.
\subsection{Invariant Differential Operators} 
\noindent
Invariant differential operators have always been a central subject of investigation  by Helgason. We find it very informative to read his first paper on the subject \cite{He59}. In retrospect, this shines a beacon to follow through much of his later work on this subject. Here we find a lucid introduction to differential operators on manifolds and the geometry of  homogeneous spaces, reminiscient of the style to appear in his famous book \cite{He62}. Specializing to a reductive homogeneous space, he begins the study of  $\D (G/H)$, those differential operators that commute with the action of the group of isometries. The investigation of this algebra of operators will occupy him through many years.  What is the relationship of $\D (G/H)$ to $\D (G)$ and what is the relationship of $\D (G/H)$ to the center of the universal enveloping algebra? Harish-Chandra had just described his isomorphism of the center of the universal enveloping algebra with the Weyl invariants in the symmetric algebra of a Cartan subalgebra, so Helgason introduces this to give an alternative description of  $\D (G/H)$.  But the goal is always  to understand analysis on the objects, so he investigates several problems, variations of which will weave throughout his research.

For symmetric spaces $X=G/K$ the algebra $\D (X)$ was known to be commutative, and Godement had formulated the notion of harmonic function in this case obtaining a mean value characterization. Harmonic functions being joint eigenfunctions of $\D (X)$ for the eigenvalue zero, one could consider eigenfunctions for other eigenvalues. Indeed, Helgason shows that the zonal spherical functions are also eigenfunctions for the mean value operator.  When $X$ is a two-point homogeneous space, and with \'Asgeirsson's result on mean value properties for solutions of the ultrahyperbolic Laplacian in Euclidean space in mind, Helgason formulates and proves an extension of it to these spaces. Here $\D (X)$ has a single generator, the Laplacian, for which he constructs geometrically a fundamental solution, thereby allowing him to study the inhomogeneous problem for the Laplacian. This paper contains still more. In many ways the two-point homogeneous spaces are ideal generalizations of Euclidean spaces so following F. John \cite{John} he is able to define a Radon like transform on the constant curvature ones and identify an inversion operator. Leaving the Riemannian case, Helgason considers harmonic Lorentz space $G/H$. He shows $\D (G/H)$ is generated by the natural second order operator; he obtains a mean value theorem for suitable solutions of the generator and an explicit inverse for the mean value operator. Finally, he examines the wave equation on harmonic Lorentz spaces and shows the failure of Huygens principle in the non-flat case.

Building on these results he subsequently examines the question of existence of fundamental solutions more generally. He solves this problem for symmetric spaces as he shows that every $D\in\D (X)$ has a fundamental solution, \cite[Thm. 4.2]{He64}. Thus, there exists a
distribution $T\in C_c^\infty (X)^\prime$ such that $DT=\delta_{x_o}$. Convolution then
provides a method to solve the inhomogeneous problem, namely, if $f\in C_c^\infty (X)$ then there exists $u\in C^\infty (X)$ such that
$Du=f$.  Those
results had been announced in  \cite{He63c}. The existence of the fundamental solution uses the deep results of Harish-Chandra on the aforementioned isomorphism as well as classic results of H\"ormander on constant coefficient operators. It is an excellent example of the combination of the classical theory with the semisimple theory. Here is a sketch of his approach.

In his classic paper \cite{HC58} on zonal spherical functions, Harish-Chandra introduced several important concepts to handle harmonic analysis. One was the appropriate notion of a Schwartz-type space of $K$ bi-invariant functions, there denoted $I(G)$. $I(G)$ with the appropriate topology is a 
Fr\'echet space, and having $C^\infty_c(X)^K$ as a dense convolution subalgebra. 

Another notion from \cite{HC58} is the Abel transform
\[ F_f(a)=a^\rho \int_N f (an)\, dn\, .\]
Today this is also called the $\rho$-twisted Radon transform
and denoted $\cR_\rho $. Eventually Harish-Chandra showed that this gives a topological isomorphism of
$I (G)$ onto $\cS (A)^W$, the Weyl group invariants in the Schwartz space on the Euclidean space $A$. Furthermore,
the Harish-Chandra isomorphism $\gamma : \D (X)\to \D (A)$ interacts compatibly in that
\[\cR_\rho (Df)=\gamma (D)\cR_\rho (f)\, .\]
One can restate this by saying that the Abel transform turns invariant differential equations
on $X$ into constant coefficient differential equations on $A\simeq \fa\simeq \R^{\rank X}$.
It follows then that $\cR_\rho^t : \cS^\prime (A)^W\to I(G)^\prime$ is also an
isomorphism. This can then be used to pull back the fundamental solution for $\gamma (D)$ to a fundamental solution for $D$.

The article \cite{He64} continues the line of investigation from \cite{He59} into the structure of $\D (X)$. 
If we denote by $U(\fg)$ the universal enveloping algebra of $\fg^\C$, then
$U(\fg)$ is isomorphic to $\D (G)$. Let $\Z (G)$ be the center of $\D (G)$. This is the algebra
of bi-invariant differential operators on $G$.
The algebra of invariant differential operators on $X$ is isomorphic to
$\D (G)^K /\D (G)^K\cap \D (G)\fk$ and therefore contains  $\Z (G)$ as an Abelian
subalgebra. 

Let $\fh$ be a Cartan subalgebra in $\fg$ extending $\fa$ and denote by
$W_\fh$ its Weyl group.   The subgroup $W_\fh (\fa)=\{w\in W_\fh\mid w(\fa)=\fa\}$ induces the little Weyl group $W$ by restriction. It follows that the restriction $p\mapsto p|_\fa$ maps
$S(\fh )^{W_\fh}$ into $S(\fa )^W$. Now $\Z (G)\simeq S(\fg^\C)^G\simeq S(\fh)^{W_\fh}$, and
$\D (X)\simeq S(\fa)^W\simeq S(\fs)^K$, $\fs$ a Cartan complement of $\fk$. The structure of these various incarnations is given in cf. \cite[Prop. 7.4]{He64} and \cite[Prop 3.1]{He92}. See also the announcements in \cite{He62a,He63c}:
\begin{theorem} The following are equivalent:
\begin{enumerate}
\item $\D (X)=\Z (G)$.
\item $S (\fh^\C)^{W_\fh}|_{\fa}=S(\fa^\C)^W$.
\item $S(\fg)^G|_{\fs}=S(\fs )^K$.
\end{enumerate}
\end{theorem}

A detailed inspection showed that (2) was always true for the classical
symmetric spaces but fails for some of the exceptional symmetric spaces. Those ideas
played an important  role in \cite{OW11} as similar restriction questions were considered
for sequences of symmetric spaces of increasing dimension.

The final answer, prompted by a question from G. Shimura, is  \cite{He92}:
\begin{theorem} Assume that $X$ is irreducible. Then $\Z (G)=\D(X)$ if and only if
$X$ is not one of the following spaces $E_6/\SO(10)\T$, $E_6/F_4$, $E_7/E_6\T$ or
$E_8/E_7\SU(2)$. Moreover, for any irreducible $X$ any $D\in\D (X)$ is a quotient of elements of $\Z(G)$.
\end{theorem}

\subsection{The Radon Transform on $\R^n$}
\noindent
The Radon transform as introduced by J. Radon in 1917 \cite{Ra17} \cite{RaGes} associates to a suitable
function $f:\R^2\to \C$ its integrals over affine lines $L\subset \R^2$
\[\cR (f)(L)=\wf (L):=\int_{x\in L} f(x)\, dx\]
for which he derived an inversion formula. This ground breaking article appeared in a not easily available journal (one can find the reprinted
article in \cite{He80}), and consequently was not well known. Nevertheless,  its true worth is easily determined by the many generalizations of it  that have been made in geometric analysis and representation theory, some already pointed out in Radon's original article. 

An important milestone in the development of the theory was F. John's book \cite{John}. Later,
the application of integration over affine lines in three dimensions played an  important role in the three dimensional X-ray transform. We refer to \cite{E03,GGG00,He80,He11,N01} for information about the history and the many applications of the Radon transform and its
descendants. 

Helgason first displayed his interest in the Radon transform in that basic paper \cite{He59}. There he considers a transform associated to totally geodesic submanifolds in a space of constant curvature and produces an inversion formula. To use it as a tool for analysis one needs to determine if there is  injectivity on some space of rapidly decreasing functions and compatibility with invariant differential operators, just as Harish-Chandra had done for the map $F_f$. In \cite{He65} Helgason starts on his long road to answering such questions, and, in the process recognizing the underlying structure as incidence geometry, he is able to describe a vast generalization. 

As he had previously considered two-point homogeneous spaces he starts there, but to this he extends Radon's case to affine p-planes in Euclidean space. We summarize the results in the important article
\cite{He65}.

Denote by  $\rH (p,n)$  the space of $p$-dimensional affine subspaces of $\R^n$. 
Let $f\in C_c^\infty (\R^n)$ and $\xi\in H(p,n)$. Define    
\[\cR(f)(\xi) = \wf (\xi ):=\int_{x\in \xi} f(x)\, d_\xi x\]
where the measure $d_\xi x$ is determined in the following way. The connected
Euclidean motion group $\rE(n)=\SO (n)\ltimes \R^n$ acts transitively on both $\R^n$ and $\rH (p,n)$.
Take basepoints  $x_o=0\in\R^n$ and  $\xi_o=\{(x_1,\ldots ,x_p, 0, \ldots ,0)\}\in\rH (p,n)$ and take
 $d_{\xi_o}x$   Lebesgue measure on
$\xi_o$. For $\xi\in \rH (p,n)$ choose  $g\in \rE (n)$
such that $\xi =g\cdot \xi_o$. Then $d_\xi x =g^*d_{\xi_o}x$ or
\[\int_\xi f (x) d_\xi x =\int_{\xi _o} f(g\cdot x)\, dx\, .\]

For $x\in \R^n$ the set $\vx:=\{\xi \in H(p,n)\mid x\in \xi\}$ is compact, in
fact isomorphic to the Grassmanian $\rG (p,n)=\SO (n)/\rS(\rO (p)\times \rO (n-p))$ of all
$p$-dimensional subspaces of $\R^n$.  Therefore each of these carries a unique
$\SO (n)$-invariant probability measure $d_x\xi$ which provides the \textit{dual Radon transform}. 
Let $\varphi \in C_c(\Xi)$  and define 
\[ \vv (x)=\int_{\vx} \varphi (\xi )\, d_x\xi\, .\]
We  have the Parseval type relationship
\[\int_\Xi \wf (\xi )\varphi (\xi)\, d\xi =\int_{\R^n} f (x)\vv  (x)\, dx\, \]
and both the Radon transform and its dual are $\rE(n)$ intertwining operators.

If $p=n-1$ every hyperplane is of the form 
$\xi =\xi (u,p)=\{x\in\R^n\mid \ip{x}{u}=p\}$ and $\xi (u,p)=\xi (v,q)$ if and only if $(u,p)=\pm (v,q)$. Thus
$\rH (p,n)\simeq \rS^{n-1}\times_{\Z_2}\R$. We now
we have the hyperplane Radon transform considered in \cite{He65}. This case had been considered by F. John \cite{John} and he  proved  the following inversion
formulas for suitable functions $f$:
\begin{eqnarray*}
f(x)&=&\frac{1}{2}\frac{1}{(2\pi i)^{n-1}}\Delta_x^{\frac{n-1}{2}}\int_{\rS^{n-1}}\wf 
(u, \ip{u}{x})\, du\, ,\quad \text{ $n$ odd}\\
f(x)&=& \frac{1}{(2\pi i)^{n}}\Delta_x^{\frac{n-2}{2}}\int_{\rS^{n-1}}
\int_{\R} \frac{\partial_p\wf 
(u, p)}{p-\ip{u}{x}}\, dp du\, ,\quad\text{ $n$ even}\, .
\end{eqnarray*}

The difference between the even and odd dimensions is significant, for in odd dimensions inversion
is given by a local operator, but not in even dimension. This is fundamental in \textit{Huygens' principle} for the wave equation to be discussed subsequently. 

For Helgason the problem is to show the existence of suitable function spaces on which these transforms are injective and to show they are compatible with the  $\rE(n)$ invariant differential operators. One shows that the Radon transform extends to the Schwartz space $\cS (\R^n)$ of
rapidly decreasing functions on $\R^n$ and it maps that space into a suitably defined Schwartz space 
$\cS (\Xi)$ on $\Xi$.  Denote by $\D (\R^n)$, respectively $\D (\Xi)$, the algebra of $\rE (n)$-invariant differential
operators on $\R^n$, respectively $\Xi$ . Furthermore, define a differential operator $\Box $ on $\Xi$ by
$\Box f(u,r) =\partial^2_r f(u,r)$.

However a new feature arises whose existence suggests future difficulties in generalizations.
Let
\[\cS^*(\R^n)=\{f\in \cS (\R^n)\mid \int_{\R^n}f(x)p(x)\, dx=0
\text{ for all polynomials } p(x)\}\]
and 
\[\cS^* (\Xi )=\{\varphi \in \cS (\Xi)\mid \int_{\R }\varphi (u,r) q(r)\, dr=0\text{ for all polynomials } q(r)\}\, .\]
Finally, let  $\cS_H (\Xi )$ be the space of rapidly decreasing function on $\Xi$ such that for each $k\in \Z^+$ the
integral $\int \varphi (u,r)r^k\, dr$ can be written as a homogeneous polynomial in $u$ of degree $k$. 
Then we have the basic theorem for this transform and its dual:
\begin{theorem}\label{th-Rad} \cite{He65} The following hold:
\begin{enumerate}
\item $\D (\R^n)=\C [\Delta]$ and $\D (\Xi) = \C [\Box ]$.
\item $\widehat{\Delta f}=\Box \wf$.
\item The Radon transform is a bijection of $\cS (\R^n)$ onto $\cS_H (\Xi)$  and the dual transform
is a bijection $\cS_H (\Xi)$ onto $\cS (\R^n)$.
\item 
The Radon transform is a bijection of $\cS^*(\R ^n)$ onto $\cS^* (\Xi)$ and the dual transform
is a bijection $\cS ^*(\Xi)$ onto $\cS^*(\R^n)$.
\item  Let $f\in \cS (\R^n)$ and $\varphi \in \cS^* (\Xi)$. If $n$ is odd then
\[f =c\, \Delta^{(n-1)/2} (\wf)^\vee \text{ and }\varphi=c\, \Box^{(n-1)/2}(\vv)^\wedge\]
for some constant independent of $f$ and $\varphi$.
\item  Let $f\in \cS (\R^n)$ and $\varphi \in \cS^* (\Xi)$. If $n$ is even then
\[f=c_1\, J_1 (\wf)^\vee \text{ and }\varphi=c_2 \, J_2 (\vv)^\wedge \]
where the operators $J_1$ and $J_2$ are given by analytic continuation
\[J_1: f(x)\mapsto \text{an.cont}\vert_{\alpha=1-2n}\int_{\R^n} f (y)
\|x-y\|^\alpha\, dy\]
and
\[J_2 : \varphi \mapsto  \text{an.cont}\vert_{\beta=-n}\int_{\R} \varphi (u,r)
\|s-r\|^\beta\, dr\]
and $c_1$ and $c_2$ are constants independent of $f$ and $g$.
\end{enumerate}
\end{theorem}

In \cite{He80} it was shown that the map
\[f\mapsto \Box^{(n-1)/4}\, \wf\]
extends to an isometry of $L^2(\R^n)$ onto $L^2(\Xi)$ .

 Needed for the proof of the theorem is one of his fundamental contributions to the subject in the following support theorem in \cite{He65}. An important generalization of this theorem will be crucial for his later work on solvability of invariant differential operators on symmetric spaces. 
\begin{theorem}[Thm 2.1 in \cite{He65}] Let $f\in C^\infty (\R^n)$ satisfy the following conditions:
\begin{enumerate}
\item For each integer $x\mapsto \|x\|^k|f(x)|$ is bounded.
\item There exists a constant $A>0$ such that $\wf (\xi)=0$ for $d(0,\xi )>A$.
\end{enumerate}

Then $f(x)=0$ for $\|x\|>A$.
\end{theorem}

An important technique in the theory of the Radon transform, which also plays an important role
in the proof of Theorem \ref{th-Rad}, uses the \textit{Fourier slice formula}:
Let $r>0$ and $u\in \rS^{n-1}$ then
\begin{equation}\label{eq-Fslides}
\cF (f)(ru)=c \int_{\R}\wf (u,s) e^{-isr}\, ds\, .
\end{equation}

So that if $f$ is supported in a closed ball $B_r^n (0)$  in $\R^n$ of radius $r$ centered at the origin, then
by the classical Paley-Wiener theorem for $\R^n$ the function
\[r\mapsto \cF (f)(ru)\]
extends to a holomorphic function on $\C$ such that
\[\sup_{z\in \C}(1+|z|)^n e^{-r|\Im z|}|\cF( f)(zu)<\infty\, \]
Let $C_{r,H}^\infty (\Xi)$ be the space of $\varphi\in \cS_H (\Xi)$ such that
$p\mapsto \varphi (u,p)$ vanishes for $p>r$.
Then the Classical Paley-Wiener theorem combined with (\ref{eq-Fslides}) shows
that the Radon transform is a bijection  $C^\infty_r(\R^n)\simeq C^\infty_{r,H}(\Xi)$, \cite[Cor. 4.3]{He65}.
(\ref{eq-Fslides})  also played an important
role in Helgason's introduction of the Fourier transform on Riemannian symmetric spaces
of the noncompact type.

\subsection{The Double Fibration Transform}
\noindent
The Radon transform on $\R^n$ and the dual transform are examples of
\textit{the double fibration transform}  introduced in \cite{He66b,He70}.
Recall that both $\R^n$ and $\rH (p,n)$ are
homogeneous spaces for the group $G=\rE (n)$.
Let $K=\SO (n)$, $L=\rS (\rO (p)\times \rO( n-p))$ 
and
$N=\{(x_1,\ldots ,x_p,0,\ldots 0)\mid x_j\in\R\}\simeq \R^p$
and $H=L\ltimes N$. Then $\R^n\simeq G/K$, $\rH (p,n)\simeq G/H$ and $L=K\cap H$.
Hence we have the double fibration

\begin{equation}\label{eq:dFibDi}
 \xymatrix{&G/L \ar[dl]_\pi \ar@{->}[dr]^p&\\
X=G/K& &\Xi=G/H } 
\end{equation}
\medskip

\noindent
where $\pi$ and $p$ are the natural projections. If $\xi=a\cdot \xi_o\in \Xi$
and $x=b\cdot x_o\in X$ then
\begin{equation}\label{eq:RtrDfib}
\wf (\xi )=\int_{H/L}f (a h\cdot x_o)\, d_{H/L}(hL)
\end{equation}
and 
\begin{equation}\label{eq:DRtrDfib}
\vv (x)=\int_{K/L} \varphi (bk\cdot \xi_o)\, d_{K/L} (kL)
\end{equation}
for suitable normalized invariant measures on $H/L\simeq N$ and $K/L$.

More generally, using Chern's formulation of integral geometry on homogeneous spaces as \textit{incidence geometry} \cite{Chern}, Helgason introduced
the following  \textit{double fibration transform}.
Let $G$ be a locally compact Hausdorff topological group  and $K,H$ two closed subgroups giving the double fibration in \ref{eq:dFibDi}.
We will assume that $G$, $K$, $H$ and $L:=K\cap H$ are all unimodular. Therefore each
of the spaces $X=G/K$, $\Xi=G/H$, $G/L$, $K/L$ and $H/L$ carry an invariant
measure.  

We set $x_o=eK$ and $\xi_o=eH$. Let $x=aK\in X$ and
$\xi=bH\in\Xi$. We say that $x$ and $\xi$ are incident if $aK\cap bH\not=\emptyset$.
For $x\in X$ and $\xi \in \Xi$ we set
\[\hat{x}=\{\eta \in \Xi\mid x \text{ and } \xi \text{ are incident }\}\]
and similarly
\[\xi^\vee =\{x \in X\mid \xi \text{ and } x \text{ are incident }\}\, .\]

Assume that if
 $ a\in K$ and $aH\subset HK$ then
$a\in H$ and similarly, if $b\in H$ and $bK\subset KH$ then $b\in K$. Thus we can
view the points in $\Xi$ as subsets of $X$, and similarly points in $X$ are
subsets of $\Xi$. Then $x^\vee$ is the set of all $\xi$ such that
$x\in \xi$ and $\hat\xi$ is the set of points $x\in X$ such that $x\in \xi$. We also have
\[\hat{x}=p (\pi^{-1}(x)) = aK\cdot \xi_0\simeq H/L
\text{ and }
\xi^\vee =\pi (p^{-1}(\xi )) =b H\cdot x_o \simeq K/L.\] 

Under these conditions the \textit{Radon transform} (\ref{eq:RtrDfib}) and
its dual (\ref{eq:DRtrDfib}) are well defined at least for compactly supported functions. Moreover,
for a suitable normalization of the measures we have 
\[\int_{\Xi}\hat f (\xi )\varphi (\xi )\, d\xi =\int_X f(x)  \varphi^\vee (x)\, dx\, .\]

Helgason \cite[p.39]{He66b} and \cite[p.147]{He84} proposed the following problems for these transforms $f \to \hat{f}$,
$\varphi \to \varphi^\vee $:

\begin{enumerate}

\item Identify function spaces on $X$ and $\Xi$  related by the integral transforms
$f\mapsto \wf$ and $\varphi \mapsto \vv$.

\item Relate the functions $f$ and $\wf^\vee$ on $X$, and similarly
$\varphi$ and $(\vv)^\wedge$ on $\Xi$, including an inversion formula, if possible.

\item Injectivity of the transforms and description of the image.

\item Support theorems.

\item For $G$ a Lie group, with $\D (X)$, resp. $\D (\Xi)$, the 
algebra of invariant differential operators on $X$, resp. $\Xi$. Do there exist 
maps $D\mapsto \widehat{D}$ and $E\mapsto E^\vee$ such that
\[(Df)^\wedge = \widehat{D}\wf\text{ and } (E\varphi )^\vee=
E^\vee \varphi^\vee\, .\]
\end{enumerate}

 There are several examples where the double fibration transform serves as a guide, e.g. the Funk transform on the sphere $\rS^n$, see \cite{F16} for
the case $n=2$, and more generally \cite{Ru02}; and  the geodesic $X$-ray transform on compact symmetric spaces,
see \cite{He07,R04}. Other uses of the approach can be found in \cite{K11}.
We refer the reader to \cite{E03} and \cite {He11} for more examples.

\subsection{Fourier analysis on $X = G/K$}
{}From now on $G$ will stand for a non-compact connected semisimple Lie group with finite center and $K$ a maximal compact subgroup. We take an Iwasawa decomposition $G=KAN$ and use standard notation for projections on to the $K$ and $A$ component. Set $X=G/K$
as before
and denote by $x_o$ the base point $eK$. Given Helgason's classic presentation of the structure of symmetric spaces \cite{He62} there is no good reason for us to repeat it here, so we use it freely and we encourage those readers new to the subject to learn it there.

In this section we introduce Helgason's version of the Fourier transform on $X$, see \cite{He65a,He68,He70}.
At first we follow the exposition in \cite{OS08} which is based more
on representation theory, i.e. \`a la von Neumann and Harish-Chandra, rather than geometry as did Helgason. For additional information see the more modern representation theory approach of \cite{OS08}, although we caution the reader
that in some places notation and definitions differ. 
 
The regular action of $G$ on $L^2(X)$ is 
$\ell_{g}f(y)=f(g^{-1}\cdot y)$, $g\in G$ and $y\in X$. 
For an irreducible unitary representation $(\pi , V_\pi)$ of $G$ and $f\in L^1(X)$ set
\[\pi (f)=\int_{G} f(g)\pi (g)\, dg\, .\]
Here we have pulled back $f$ to a right $K$-invariant function on
$G$. If $\pi (f)\not= 0$ then $V_\pi^K=
\{v\in V_\pi \mid (\forall k\in K)\, \pi (k)v=v\}$ is nonzero. Furthermore,
as $(G,K)$ is a Gelfand pair we have $\dim V_\pi^K =1$, in which case $(\pi, V_\pi)$ is called spherical.

Fix a unit vector $e_\pi\in V_\pi^K$. Then 
$\Tr (\pi (f)) =(\pi (f)e_\pi,e_\pi)$ and $\|\pi (f)\|_{\text{HS}}=\|\pi (f)e_\pi\|$.
Note that both $(\pi (f)e_\pi, e_\pi)$ and $\|\pi (f)e_\pi\|$ are
independent of the choice of $e_\pi$.  Let
$\wG$ be the set of equivalence classes of irreducible unitary spherical
representations of $G$. Then as $G$ is a type one group, there exists a
measure $\hm$ on $\wG$ such that
\begin{equation}\label{eq:dirInt}
f(g\cdot x_o)=\int_{\wG} (\pi (f)e_\pi ,\pi(g)e_\pi)d\hm (\pi)
\text{ and } 
\|f\|_2^2=\int_{\wG}\|\pi (f)e_\pi \|^2_{\text{HS}}\, d\hm (\pi )\, .\end{equation}

Harish-Chandra, see  \cite{HC54,HC57,HC58,HC66}, determined the representations that occur in
the support of the measure in the decomposition (\ref{eq:dirInt}),  as well as an explicit formula for the Plancherel measure for
the spherical Fourier transform defined by him.

Helgason's formulation is motivated by ``plane waves". First we fix parameters. Let $(\lambda ,b)\in \fa_\C^*\times K/M$ and define an ``exponential function''
$e_{\lambda, b} : X\to \C$ by 
\[e_{\lambda, b}(x)=e_b (x)^{\lambda-\rho}\, ,\]
where $e_b(x)=a(x^{-1}b)$ from the Iwasawa decomposition. 
Let
$\cH_\lambda =L^2(K/M)$ with action
\[\pi_\lambda (g)f (b)=e_{\lambda, b}(g\cdot x_o)  f(g^{-1} \cdot b)\, .\]
It is easy to see that $\pi_\lambda$ is a representation with a
$K$-fixed vector $e_\lambda (b)=1$ for all $b\in K/M$; there is a $G$-invariant pairing
\begin{equation}\label{eq:pair}
\cH_{\lambda}\times \cH_{-\bar \lambda}\to \C\, ,\quad
\ip{f}{g}:=\int_{K/M} f(b)\overline{g(b)}\, db\, ;
\end{equation}
and it is unitary if and only if $\lambda \in i\fa^*$ and irreducible
for almost all $\lambda$ \cite{Kost,He76}.
With $\wtf_\lambda :=\pi_{\lambda}(f)e_{\lambda}$ we have
\begin{equation}\label{eq:FT}
\wtf_\lambda (b) = \pi_\lambda (f)e_\lambda (b)=
\int_G f(g) \pi_\lambda (g)e_\lambda (b)  \, dg
=\int_X f(x) e_{\lambda ,b} (x)\, dx  \, .
\end{equation}
Then $\wtf (\lambda ,b):=\wtf_{i\lambda} (b)$ is the Helgason
Fourier transform on $X$, see \cite[Thm 2.2]{He65a}.

Recall the little Weyl group $W$. The representation $\pi_{w\lambda}$
is known to be equivalent with $\pi_\lambda$ for almost all
$\lambda\in \fa^*_\C$. Hence for such $\lambda$ there exists an intertwining operator
$$\cA(w,\lambda)\colon \cH_\lambda \to \cH_{w\lambda}\, .$$
The operator is unique, up to scalar multiples,
by Schur's lemma. We normalize it so that
$\cA (w,\lambda )e_\lambda =e_{w\lambda}$. The family $\{\cA (w,\lambda)\}$  depends meromorphically on $\lambda$ and $\cA (w,\lambda)$
is unitary for $\lambda\in i\fa^*$. Our normalization implies that
\begin{equation}\label{eq:intTwRel}
\cA (w,\lambda )\wtf_\lambda =\wtf_{w\lambda}\, .
\end{equation}

We can now formulate the Plancherel Theorem for the Fourier transform in
the following way, see \cite[Thm 2.2]{He65a} and also \cite[p. 118]{He70}.

First let 
\[\bc (\lambda )=\int_{\bar N} a(\bar n)^{-\lambda -\rho}\, d\bar n \]
be the Harish-Chandra $c$-function, $\lambda$ in a positive chamber. 
The Gindikin- Karpelevich formula for the $c$-function \cite{GK62} gives a
meromorphic extension of $\bc$ to all of $\fa_\C^*$. Moreover $\bc$ is regular
and of polynomial growth on $i\fa^*$.

To simplify the notation let $d\mu (\lambda ,kM)$ be the measure
$(\# W |\bc (\lambda )|)^{-1}\,d\lambda  d(kM) $ on $i\fa^*\times K/M$.:
\begin{theorem}[\cite{He65a}] \label{th:Pla}
The Fourier transform establishes a
unitary isomorphism
\[L^2(X)\simeq \int^{\oplus}_{i\fa^*/W}(\pi_\lambda ,\cH_\lambda )\,
\frac{d\lambda}{|\bc (\lambda )|^{2}}\, .\]
Furthermore, for $f\in C_c^\infty (X)$ we have
\[f(x)=\int_{i\fa^*\times {K/M}} \wtf_\lambda (b) e_{-\lambda,b}(x)\, d\mu (\lambda ,b) \, .\]
\end{theorem}

Said more explicitly, the Fourier transform extends to
a unitary isomorphism
\begin{eqnarray*} L^2(X)&\to& L^2\left(i\fa^*,d\mu, L^2(K/M)\right)^W\\
&=&\left\{\left.\varphi \in L^2\left(i\fa^*,d\mu ,L^2 (K/M)\right)\,\right|\, 
(\forall w\in W) \cA (w,\lambda )\varphi (\lambda) = \varphi (w\lambda) \right\}\, .
\end{eqnarray*}
To connect it with Harish-Chandra's spherical transform notice that if $f$ is left $K$-invariant, then
$b\mapsto \wtf_\lambda (b)=
\wtf (\lambda)$ is
independent of $b$ and  the integral \ref{eq:FT} can be written
as
\begin{equation}\label{eq:SphFT}
\wtf (\lambda )=\int_X\int_{K/M} e_{\lambda ,b}(x)\, db\, dx
=\int_X f(x)\varphi_\lambda (x)\, dx
\end{equation}
where $\varphi_\lambda$ is the spherical function 
\[\varphi_\lambda (x)=\int_K a(g^{-1}k)^{\lambda -\rho} \, dk\, .\]
Then (\ref{eq:SphFT}) is exactly the Harish-Chandra spherical
Fourier transform \cite{HC58} and the proof of Theorem \ref{th:Pla} can
be reduced to that formulation. 

Since $\varphi_\lambda =\varphi_\mu$ if and only if that there exists
$w\in W$ such that $w\lambda =\mu$, the spherical Fourier transform
$\wtf (\lambda )$ is $W$ invariant. The Plancherel Theorem 
reduces to
\begin{theorem}
The spherical Fourier transform sets up an unitary
isomorphism
\[L^2(X)^K\simeq L^2\left(i\fa^*/W , \frac{d\lambda}{|\bc (\lambda )|^2}\right)\, .\]
If $f\in C_c(X)^K$ then
\[f(x)=\frac{1}{\# W}\int_{i\fa^*} \wtf (\lambda )\varphi_{-\lambda}(x)
\frac{d\lambda}{|\bc (\lambda )|^2}\, .\]
\end{theorem}

A very related result is the Paley-Wiener theorem
which describes the image of the smooth compactly supported functions by the Helgason Fourier Transform. For $K$-invariant
functions  in \cite{He66} Helgason formulated the problem and solved it modulo an interchange of a specific
integral and sum. The justification for the interchange was provided in \cite{G71};  a new proof was given
in \cite[Ch.II Thm. 2.4]{He70}. The Paley-Wiener theorem for functions in $C_c^\infty (X)$ was
announced in \cite{He73a} and a complete proof was given in \cite[Thm. 8.3]{He73}. Later,   Torasso \cite{T77} produced another proof, and Dadok  \cite{D79}  generalized it to
distributions on $X$. 

There are many applications of the Paley-Wiener Theorem and the ingredients of its proof. For example an alternative approach to the inversion formula can be obtained \cite{Ros}. The Paley-Wiener theorem was used in \cite{He73} in the  proof of surjectivity discussed in the next section, and in \cite{He76} to prove the necessary and sufficient condition for the bijectivity
of the Poisson transform for $K$-finite functions on $K/M$ to be discussed subsequently. The Paley-Wiener theorem plays an important role in the
study of the wave equation on $X$ as will be discussed later.

For the group $G$, an analogous theorem, although much more complicated in statement and proof, was finally obtained by Arthur
\cite{A83}, see also \cite{CD84,CD90,vBS05}. In \cite{De} the result was extended to non $K$-finite functions. The equivalence of the apparently different formulations of the characterization can be found in \cite{BanSou}. For semisimple symmetric
spaces $G/H$ it was done by van den Ban and Schlichtkrull
\cite{vBS06}. The local Paley-Wiener theorem for compact groups was
derived by 
Helgason's former student F. Gonzalez in \cite{G01} and then for all compact
symmetric spaces in
\cite{BOP95a,C06,OS08,OS10,OS11}. 

\subsection{Solvability for $D\in\D (X)$}
We come to one of Helgason's major results: a resolution of the solvability problem for $D\in\D (X)$. We have seen the existence of a fundamental solution allows one to solve the inhomogeneous equation: given $f\in C_c^\infty (X)$ does there exists $u\in C^\infty (X)$  with $Du=f?$ But what if $f\in C^\infty (X)$? This is much more difficult. Given Helgason's approach outlined earlier it is natural that once again he needs a Radon-type transform but more general than for $K$ bi-invariant functions.

The Radon transform on symmetric spaces of the noncompact type is, 
as mentioned in the earlier section, an example of the double fibration transform and probably one
of the motivating examples for S. Helgason to introduce this general framework.
Here the
double fibration is given by

\begin{equation}\label{eq:dFibDi2}
 \xymatrix{&G/M \ar[dl]_\pi \ar@{->}[dr]^p&\\
X=G/K& &\Xi=G/MN } 
\end{equation}
\medskip

\noindent
and the corresponding transforms are for compactly supported functions:
\[\wf (g\cdot \xi_o)=\int_N f(gn\cdot x_o)\, dn
\text{ and } \vv (g\cdot x_o)=\int_K \varphi (gk\cdot \xi_o)\, dk\, .\]

As mentioned before, in the $K$ bi-iinvariant setting this type of Radon transform
had already appeared (with an extra factor $a^\rho$) in the work of Harish-Chandra 
\cite{HC58} via the map $f\mapsto F_f$. It also appeared in the fundamental work
by Gelfand and Graev \cite{GG59, GG62} where they introduced the ``horospherical method''.

In this section we introduce the Radon transform on $X$ and discuss some of
its properties. It should be noted that Helgason introduced the Radon
transform in \cite{He63a,He63b} but the Fourier transform only appeared later in \cite{He65a}, see also \cite{He66b}.

We have seen that the Fourier transform on $X$ gives a unitary
isomorphism
\[L^2(X)\simeq \int_{\fa^+}^\oplus (\pi_\lambda ,\cH_\lambda )\, \frac{d\lambda}{|\bc (\lambda )|^2}\]
whereas the Fourier transform in the $A$-variable gives a unitary isomorphism
\[L^2(\Xi )\simeq \int^{\oplus}_{i\fa} (\pi_\lambda ,\cH_\lambda )\, d\lambda\, .\]

As the representations $\pi_{\lambda}$ and $\pi_{w\lambda}$, $w\in W$, are
equivalent this has the equivalent formulation
\[L^2(\Xi) \simeq (\# W ) L^2(X)\, .\]

In hindsight we could 
construct an intertwining operator from the following sequence of
maps
\[L^2(X)\to L^2\left(K/M\times i\fa^*,\frac{d\lambda}{|\bc (\lambda )|^2}\right)
\to L^2(K/M\times i\fa^*, d\lambda )\to L^2(\Xi )\]
obtained with $b=k\cdot b_o$ from the sequence:
\[f\mapsto \wtf_\lambda (b) \mapsto
\frac{1}{\bc (\lambda )}\wtf (\lambda,b)
\mapsto \cF_A^{-1} \left(\frac{1}{\bc (\cdot  )}\wtf(\cdot,b)\right)(a)
=:\Lambda (f)(ka\cdot \xi_o )\, .\]

This idea plays a role in the inversion of the Radon transform, but instead we start with the
Fourier transform on $X$ given by (\ref{eq:FT}). Then using $b=k\cdot b_o$ we have
\begin{eqnarray*}
\wtf (\lambda,b )&=&\int_X f(x)e_{\lambda ,b}(x)\, dx\\
&=&\int_X f(g\cdot x_o) a(g^{-1}l)^{\lambda - \rho}\, dg\\
&=& \int_X f (lg\cdot x_o) a(g^{-1})^{\lambda -\rho}, dg\\
&=& \int_A\int_N f(lan\cdot x_o)a^{-\lambda +\rho} \, dnda\\
&=&\cF_A ((\cdot )^{\rho }\cR(f)(l (\cdot ))(\lambda)\, .\end{eqnarray*}

Here $\cR (f)=\hat f$ is the \textit{Radon Transform} from before. 
Thus we obtain that the factorization of the
unitary $G$ map discussed above, namely the Fourier transform on $L^2(X)$ is followed by the Radon transform, which is then followed by
the Abelian Fourier transform on $A$, all this modulo the application
of the pseudo-differential operator $J$ corresponding to the Fourier multiplier 
$1/\bc (\lambda )$. Following 
\cite{He65a}  and \cite[p. 41 and p. 42]{He70} we therefore define the
operator $\Lambda$ by
\[\Lambda (f)(ka\cdot \xi_o) =a^{-\rho}J_a (a^\rho f(k a\cdot \xi_o))\, .\]

 We then get \cite[Thm. 2.1]{He65a} and  \cite{He70}:

\begin{theorem} Let $f\in C_c^\infty (X)$. Then
\[\#W \int_X |f(x)|^2\, dx =\int_{\Xi} |\Lambda \cR (f) (\xi )|^2\, d\xi\]
and $f\mapsto \frac{1}{\# W} \Lambda \cR (f)$ extends to an isometry into
$L^2(X)$. Moreover, for $f\in C_c^\infty (X)$
\[f(x)= \frac{1}{\# W} (\Lambda \Lambda^*\hat f)^\vee (x)\, .\] 
\end{theorem}
With inversion in hand, in \cite{He63b} and \cite{He73} Helgason obtains the key properties of the Radon transform needed for the analysis of invariant differential operators on $X$. First we have the compatibility with a type of Harish-Chandra isomorphism:
\begin{theorem} There exists a homorphism $\Gamma : \D (X)\to \D (\Xi)$ such
that for $f\in C_c(X)$ we have $\cR (D f)=\Gamma (D)\cR (f)$.
\end{theorem}

Then using the Paley-Wiener Theorem for the symmetric space $X$  Helgason generalizes his earlier support theorem.

\begin{theorem}[ \cite{He73}]

Let $f\in C_c^\infty (X)$ satisfy the following conditions:
\begin{enumerate}
\item There is a closed ball $V$ in $X$.
\item The Radon transform $\wf (\xi)=0$ whenever the horocycle $\xi$ is disjoint from $V$.
\end{enumerate}

Then $f(x)=0$ for $x\notin V$.
\end{theorem}

He now has all the pieces of the proof of his surjectivity result.

\begin{theorem}\cite[Thm. 8.2]{He73}    Let $D\in \D (X)$. Then
$$ D C^\infty (X) = C^\infty (X).$$
\end{theorem}

The support theorem has now been extended to noncompact reductive symmetric spaces by Kuit \cite{K11}.

\subsection{The Poisson Transform}
On a symmetric space $X$ the use of the Poisson transform has a long and rich history. But into this story fits a very precise and important contribution - the\,  \lq\lq Helgason Conjecture". In this section we recall briefly the background from Helgason's work leading to this major result.
 
Let $g\in L^2(K/M)$ and $f\in C_c^\infty (X)$. Recall from Theorem \ref{th:Pla} that the Fourier transform can be viewed as having values in $ L^2(i\fa^*,\frac{d\lambda}{\# W|\bc (\lambda )|^{2}}, L^2(K/M))^W.$ Denote the Fourier transform on $X$  by $\cF_X(f)(\lambda )=\wtf_\lambda$ and by $\cF_X^*$ its adjoint. Then we
evaluate $\cF_X^*$ as follows
\begin{eqnarray*}
\ip{\cF_X(f)}{g}&=&\ip{f}{\cF_X^*(g)}\\
&=&\int_Xf(x) \int_{i\fa^*}\left(\overline{\int_{K/M} e_{- \lambda ,b} (x)g(b)\, db}\, \right)
\frac{d\lambda}{|c(\lambda)|^2}\, dx\, .
\end{eqnarray*}
The function inside the parenthesis is the \textit{Poisson transform}
\begin{equation}\label{eq:PoissonTr}
\cP_{\lambda }(g)(x):=\int_{K/M} e_{-\lambda,b}(x) g(b)\, db.
\end{equation}
Helgason had made the basic observation that the functions $e_{\lambda, b}$ are eigenfunctions for $\D (X)$, i.e., there exists a character $\chi_\lambda : \D (X)\to \C$ such that
\[De_{\lambda ,b}= \chi_\lambda (D) e_{\lambda,b}\, .\]
Indeed, they are fundamental to the construction of the Helgason Fourier transform. Here they form the kernel of the construction of eigenfunctions.

Let 
\begin{equation}\label{eq:Elambda}
\cE_\lambda (X):=\{f\in C^\infty (X)\mid (\forall D\in\D (X))\, Df=\chi_\lambda (D)f\}\, .
\end{equation}

Since $D\in \D(X)$ is invariant the group $G$ acts on $\cE_\lambda$.  This defines a continuous representation of $G$ 
where $\cE_\lambda$ carries the topology inherited from $C^\infty (X)$.
We have $\cP_\lambda g\in \cE_\lambda$ and
$\cP_\lambda :\cH_\lambda^\infty=C^\infty (K/M)\to \cE_\lambda$ is
an intertwining operator.

In the basic paper \cite{He59} we have seen that various properties of joint solutions of operators in $\D(X)$ are obtained. In hindsight, one might speculate about eigenvalues different than $0$ for operators in $\D(X)$, and what properties the eigenspaces might have. In fact, such a question is first formulated precisely in \cite{He70} where several results are obtained. Are the eigenspaces irreducible? Do the eigenspaces have boundary values? What is the image of the Poisson transform on various function spaces?

In \cite{He70} Helgason observed that, as $b\mapsto
e_{-\lambda ,b}(x)$ is analytic, the Poisson transform extends to the dual $\cA^\prime (K/M)$ of
the space $\cA (K/M)$ of analytic functions on $K/M$. 

 Recall the Harish-Chandra $c$-function $\bc (\lambda )$ and denote by 
$\Gamma_X(\lambda )$ the denominator of  $\bc( \lambda )\bc (-\lambda )$. The Gindikin- Karpelevich formula for the $c$-function  gives an explicit formula for $\Gamma_X(\lambda )$ as a product of $\Gamma$-functions. An element $\lambda \in \fa_\C^*$ is \textit{simple} if the
Poisson transform $\cP_\lambda : C^\infty (K/M,)\to \cE_\lambda (X)$ is injective. 
 \begin{theorem} [Thm. 6.1 \cite{He76}] $\lambda$ is simple if and only if
the denominator of the Harish-Chandra
$c$-function is non-singular at $\lambda$.
\end{theorem}

This result was used by Helgason for the following criterion for irreducibility:

\begin{theorem}[Thm 9.1,  Thm. 12.1, \cite{He76}]\label{elambda} The following are equivalent:
\begin{enumerate}
\item The representation of $G$ on
$\cE_\lambda(X)$ is irreducible.
\item The principal series representation $\pi_\lambda$ is irreducible.
\item $\Gamma_X (\lambda)^{-1} \not= 0$.
\end{enumerate}
\end{theorem}

In \cite{He76} p.217 he explains in detail the relationship of this result to \cite{Kost}. With irreducibility under control, Helgason turns to the range question. In \cite{He76} for all symmetric spaces of the non-compact type, generalizing
 \cite[Thm. 3.2]{He70} for rank one spaces, he proves

\begin{theorem} 
Every $K$-finite function in $\cE_\lambda(X)$ is of the form
$\cP_\lambda (F)$ for some $K$-finite function on $K/M$.
\end{theorem}

In \cite[Ch. IV,Thm. 1.8]{He70} he examines the critical case of the Poincar\'e disk. Utilizing classical function theory on the circle he shows that eigenfunctions have boundary values in the space of analytic functionals. This, coupled with the aforementioned analytic properties of the Poisson kernel allow him to prove 
\begin{theorem} $\cE_\lambda(X)=\cP_\lambda (\cA'(\T ))$ for $\lambda\in i\fa^*$
\end{theorem}

Those results initiated intense research related to finding a suitable compactification of $X$ compatible with eigenfunctions of $\D(X)$; to hyperfunctions as a suitable class of objects on the boundary to be boundary values of eigenfunctions; to the generalization of the Frobenius regular singular point theory to encompass the operators in $\D(X)$; and finally to the analysis needed to treat the Poisson transform
and eigenfunctions on $X$. The result culminated in the impressive proof by Kashiwara, Kowata, Minemura, Okamoto,
Oshima and Tanaka \cite{KKMOOT78} that the Poisson transform is a surjective map
from the space of hyperfunctions on $K/M$ onto $\cE_\lambda (X)$, referred to as the \lq\lq Helgason Conjecture".

\subsection{Conical Distributions}
\noindent
Let $X$ be the upper halfplane $\C^+=\{z\in\C\mid \Re (z)>0\}=\SL(2,\R)/\SO (2)$.
A horocycle in $\C$ is a circle in $X$ meeting the real line tangentially or, if the point of tangency is $\infty$,  real lines parallel to the $x$-axis. It is easy to see that the horocycles are the orbits of conjugates of
the group 
\[N=\left\{\left. \begin{pmatrix} 1 & x\\ 0 & 1\end{pmatrix}
\, \right|\, x\in \R\right\}\, .\]
This leads to the definition for arbitrary  symmetric spaces of the
noncompact type:
\begin{definition} A horocycle in $X$ is an orbit of a conjugate of $N$.
\end{definition}
Denote by $\Xi$ the set of horocycles.
Using the Iwasawa decomposition it is easy to see that the horocycles are 
the subsets of $X$ of the form $gN\cdot x_o$. Thus $G$ acts transitively on
$\Xi$ and $\Xi=G/MN$. As we saw before
\begin{equation}\label{eq-L2Xi}
L^2(\Xi) \simeq \int^{\oplus}_{i\fa^*} (\pi_\lambda ,\cH_\lambda )\, d\lambda \simeq (\# W ) L^2(X)
\end{equation}
the isomorphism being
given by
\[\wphi_\lambda (g):=\int_A [a^{\rho}\varphi (ga\cdot \xi_o)]a^{-\lambda}\, da
=\cF_A ([(\cdot )^\rho \ell_{g^{-1}}\varphi ] |_A)(\lambda) \, .\]

The description of 
$L^2(\Xi) \simeq (\# W ) L^2(X)$
suggests the question of relating $K$ invariant vectors with $MN$ invariant vectors.  But, as $MN$ is noncompact,  it follows from the
theorem of Howe and Moore \cite{HM79} that the unitary
representations $\cH_\lambda$, $\lambda\in i\fa^*$ do not have any nontrivial $MN$-invariant
vectors.  But they have $MN$-fixed distribution vectors as we will explain.

Let $(\pi,V_\pi)$ be a representation of $G$ in the Fr\'echet space $V_\pi$.  Denote by $V_\pi^\infty$ the
space of smooth vectors with the usual Fr\'echet topology. The space $V_\pi^\infty$ is invariant under $G$ and we denote
the corresponding representation of $G$ by $\pi^{\infty}$. The conjugate linear dual of $V_\pi^\infty$ is denoted by
$V_\pi^{-\infty}$. The dual pairing
$V^{-\infty}_\pi \times V^{\infty}_\pi\to \C$,  is denoted $\ip{\cdot }{\cdot }$. The group $G$ acts on $V_\pi^{-\infty}$ by
\[\ip{\pi^{-\infty}(a)\Phi} {\phi}:=\ip{\Phi}{\pi^\infty (a^{-1})\phi}\, .\]

The reason to use the conjugate dual is so that for unitary representations $(\pi ,V_\pi)$ we have canonical $G$-equivariant inclusions
\[V_\pi^{\infty}\subset V_\pi  \subset V_\pi^{-\infty}\, .\] 
For the principal series representations we have  more generally by (\ref{eq:pair}) $G$-equivariant embeddings
$\cH_{\bar \lambda}\subset \cH_{-\lambda}^{-\infty}$.

Assume that there exists a nontrivial distribution vector $\Phi\in ( V_\pi^{-\infty})^{ MN}$. Then we  define
$T_\Phi : V_\phi^\infty \to C^\infty (\Xi)$ by $T_\Phi (v;g\cdot \xi_o)= \ip{\pi^{-\infty}(g)\Phi}{v}$. Similarly,
if $T: V_\pi^\infty \to C^\infty (\Xi )$ is a continuous intertwining operator we can define a $MN$-invariant
distribution vector $\Phi_T: V_\pi^\infty \to \C$ by $\ip{\Phi_T}{v}=T(v;\xi_o)$. Clearly those two maps are
inverse to each other. The decomposition of $L^2(\Xi)$ in (\ref{eq-L2Xi}) therefore suggests that for
generic $\lambda$ we should have $\dim (\cH_\lambda^{-\infty} )^{MN}=\# W$.

As second motivation for studying $MN$-invariant distribution
vectors is the following. Let $(\pi,V_\pi)$ be an irreducible unitary representation
of $G$ (or more generally an irreducible admissible representation) and let
$\Phi,\Psi \in (V_\pi^{-\infty})^{MN}$. If $f\in C_c^\infty (\Xi )$ then $\pi^{-\infty}(\overline f)\Phi$ is well defined and an element in $V_\pi^\infty$. Hence $\ip{\Psi}{ \pi^{-\infty}(\overline f)\Phi}$ is a well defined
$MN$-invariant  distribution on $\Xi$ and all the invariant differential differential operators on $\Xi$ 
coming from the center of the universal enveloping algebra act on this distribution by scalars.

A final motivation for Helgason to study $MN$-invariant distribution vectors is the construction
of intertwining operators between the representations $(\pi_{\lambda},\cH_\lambda)$ and
$(\pi_{w\lambda},\cH_{w\lambda})$, $w\in W$. This is done in Section 6 in \cite{He70} but
we will not discuss this here but refer to \cite{He70} as well as \cite{S68,KS71,KS80,VW90}
for more information. 

We now recall Helgason's construction for the principal series represenations
$(\pi_\lambda ,\cH_\lambda)$. For that it is needed that $\cH_\lambda=L^2(K/M)$ is
independent of $\lambda$ and $\cH_\lambda^\infty =C^\infty (K/M)$.
Let $m^*\in N_K(\fa )$ be such that $m^*M\in W$ is the longest element.
Then the Bruhat big cell, $\bar Nm^*AMN$,  is open and dense. Define
\begin{equation}\label{eq:psi}
\psi_{ \lambda} (g)=\left\{\begin{matrix}
a^{ \lambda -\rho} &\text{if}& g= n_1 m^*a m a n_2\in  N m^*MAN\\
0  &\text{if}& \text{otherwise.}\end{matrix}\right.\, 
\end{equation}  

If $\Re \lambda >0$ then
$\psi_\lambda \in \cH_{-\bar \lambda}^{-\infty}$ is an $MN$-invariant distribution vector. Helgason then shows in Theorem 2.7 that $\lambda \mapsto \psi_\lambda
\in \cH_{-\bar\lambda}^{-\infty}$ extends to a meromorphic family of distribution
vectors on all of $\fa_\C^*$.  Similar construction works for the other $N$-orbits 
$NwMAN$, $w\in W$, leading to distribution vectors $\psi_{w, \lambda}$.

Denote by $\D( \Xi)$ the algebra of $G$-invariant differential operators on
$\Xi$. Then $H\mapsto D_H$ extends to an isomorphisms of algebras
$S (\fa)\simeq \D (\Xi)$, see \cite[Thm. 2.2]{He70}. 

\begin{definition} A distribution $\Psi$ (conjugate linear) on $G$ is conical if
\begin{enumerate}
\item $\Psi $ is $MN$-biinvariant.
\item $\Psi$ is an eigendistribution of $\D (\Xi)$.
\end{enumerate}
\end{definition}

The distribution vectors $\psi_{w, \lambda}$ then leads to conical distributions $\Psi_{w,\lambda}$ and
it is shown in \cite{He70,He76} that those distributions generate the space of conical distributions
for generic $\lambda$.

For $\lambda\in\fa_\C^*$ let $C_c^\infty(\Xi)_\lambda^\prime$ (with the relative strong topology) denote the 
joint distribution eigenspaces of $\D (\Xi)$ containing the function
$g\cdot \xi_o\mapsto a(x)^{\lambda -\rho}$. Then $G$ acts on $C_c^\infty(\Xi)_\lambda^\prime$
and according to  \cite[Ch. III, Prop. 5.2]{He70} we have:

\begin{theorem} The representation on $C_c^\infty (\Xi)_\lambda^\prime$ is irreducible if
and only if $\pi_\lambda$ is irreducible.
\end{theorem}

\subsection{The Wave Equation}
\noindent
Of the many invariant differential equations on $X$  the wave equation frequently was the focus Helgason's attention. We shall discuss some of this work, but will omit his
later work on the \textit{multitemporal} wave equation \cite{He98a,HeS99}.

Let $\Delta_{\R^n}=\sum_{j=1}^n \dfrac{\partial^2}{\partial x_j^2}$ denote the
Laplace operator on $\R^n$. The \textit{wave-equation} on $\R^n$ is the Cauchy problem
\begin{equation}\label{eq:Wave}
\Delta_{\R^n}u (x,t)=\dfrac{\partial^2}{\partial t^2} u(x,t)\, \quad u(x,0)=f(x),\, \dfrac{\partial}{\partial t}u(x,0)=g(x)
\end{equation}
where the initial values $f$ and $g$ can be from $C_c^\infty (X)$ or another ``natural''
function space. Assume that $f,g\in C_c^\infty (\R^n)$ with support contained
in a closed ball $\overline{B}_R(0)$ of radius $R>0$ and centered at zero. The solution has a finite propagation speed in the sense that $u(x,t)=0$ if $\|x\|-R\ge |t|$. The \textit{Huygens' principle} asserts that $u(x,t)=0$ for $|t|\ge \|x\|+R$. It always holds for $n>1$ and odd but fails in even dimensions.  It holds for $n=1$ if $g\in C^\infty_c(\R)$ with mean zero. 

This equation can be considered for any Riemannian or pseudo-Riemannian
manifold. In particular it is natural to consider the wave equation for Riemannian symmetric
spaces of the compact or noncompact type. Helgason was interested in the wave equation and
the Huygens' principle from early on in his mathematical career, see \cite{He64,He77,He84a,He86,He92a,He98}.  One can probably trace that
interest to his friendship with L. \'Asgeirsson, an Icelandic mathematican who studied with Courant in
G\"ottingen and had worked on the Huygens' principle on $\R^n$.

One can assume that in (\ref{eq:Wave}) we have $f=0$ and for simplicity assume that
$g$ is $K$-invariant. Then $u$ can also be taken $K$-invariant. It is also more natural to
consider the shifted wave equation
\begin{equation}\label{eq:WaveSh}
(\Delta_{X}+\|\rho \|^2) u (x,t)=\dfrac{\partial^2}{\partial t^2} u(x,t)\, \quad u(x,0)=0,\, \dfrac{\partial}{\partial t}u(x,0)=g(x)
\end{equation}
There are three main approaches to the problem. The first is to
use the Helgason Fourier transform to reduce (\ref{eq:WaveSh}) to the
differential equation
\begin{equation}
\dfrac{d^2}{dt^2}\widetilde{u}(i \lambda ,t)=- \|\lambda \|^2 \widetilde{u}(\lambda, t)\, ,\quad 
\widetilde{u}(\lambda ,0)=0 \text{ and }
\dfrac{d}{dt}\widetilde{u}(\lambda ,0)
=\widetilde{g}(\lambda )
\end{equation}
for $\lambda \in i\fa^*$. From the inversion formula we then get
\[u(x,t)=\frac{1}{\#W }\int_{i\fa^*} \widetilde{g}(\lambda )
\varphi_\lambda (x) \frac{\sin \|\lambda\|t}{\|\lambda \|}\, \frac{d\lambda}{|\bc (\lambda )|^2}\, .\]
One can then use the Paley-Wiener Theorem to shift the path of integration. Doing that
one might hit the singularity of the $\bc (\lambda)$ function. If all the root multiplicities are
even, then  $1/\bc(\lambda )\bc (-\lambda )$ is a $W$-invariant polynomial and hence corresponds to an invariant differential
operator on $X$.
 
Another possibility is to use the Radon transform and its compatibility with invariant operators
\[\cR ((\Delta +\|\rho\|^2)f)|_A
=\Delta_A \cR (f)|_A\]
then use the Helgason Fourier transform, and finally the Euclidean result on
the Huygens' principle. This was the method used in \cite{OS92}.

Finally, in \cite{He92a} Helgason showed that
\[\frac{\sin \|\lambda \|t}{\|\lambda\|}=
\int_Xe_{i\lambda,b}(x)\, d\tau_t(x)=\int_X \varphi_{-\lambda}(x)\, d\tau_t(x)\]
for certain distribution $\tau_t$ and then proving a support theorem for $\tau_t$.

The result is \cite{OS92,He92a}:
\begin{theorem} Assume that all multiplicities are even. Then Huygens's principle holds if $\rank X$ is odd. 
\end{theorem}

It was later shown in \cite{BOS95} that in general the solution has a specific exponential
decay. In \cite{BO97} it was shown, using symmetric space duality, that the
Huygens' principle holds locally for a compact symmetric spaces if and only it holds
for the noncompact dual. The compact symmetric spaces were then treated more directly
in \cite{BOP95a}.

\begin{ack} The authors want to acknowledge the work that the referee did for this paper. His thorough and conscientious report was of great value to us for the useful corrections he made and the helpful suggestions he offered. 
\end{ack}

\end{document}